\documentclass[10pt, a4paper]{amsart}
\usepackage[utf8]{inputenc}
\usepackage{tikz-cd}
\usepackage{xcolor,colortbl}
\usepackage{latexsym,amssymb}
\usepackage[english]{babel}
\usepackage{bm}
\usepackage{amsfonts, amsthm}
\usepackage{amsmath}
\usepackage{mathrsfs}
\usepackage{multirow}
\usepackage{bm}
\usepackage{pgf,tikz}
\usepackage{hyperref}
\usepackage{mathtools}
\usepackage{comment}

\newtheorem{proposition}{Proposition}[section]

\theoremstyle{definition}

\newtheorem{corollary}[proposition]{Corollary}
\newtheorem{theorem}[proposition]{Theorem}

\theoremstyle{definition}

\newtheorem{definition}[proposition]{Definition}
\newtheorem{example}[proposition]{Example}

\def\Z{\mathbb{Z}}

\numberwithin{equation}{section}

\usepackage{mathtools}

\def\Z{\mathbb{Z}}

\def\H{\mathrm{H}}

\renewcommand{\phi}{\varphi}

\begin{document}
		\title[]{Completing the Existence Problem for Integer Relative Heffter Arrays $\H_k(n;k)$}
	
	\author[L. Mella]{Lorenzo Mella}
	\address{}
	\email{lorenzo.mella@unibs.it}

\begin{abstract}
Heffter arrays, introduced by Archdeacon \cite{A}, are combinatorial structures with applications to cyclic cycle systems and biembeddings of graphs on surfaces. Costa, Morini, Pasotti and Pellegrini proposed in \cite{CMPP} the notion of relative Heffter arrays as a generalization of classical Heffter arrays, inspired by the concept of relative difference families. 

	In their work, the existence problem of integer relative Heffter arrays $\H_k(n;k)$ was solved for all $k\neq 5$, while the case $k=5$ and $n \equiv 0 \pmod{4}$ remained open, apart from two sporadic examples with $n=8,16$. In this article, we consider this open problem and we construct an $\H_5(n;5)$ for every $n\equiv 0 \pmod{4}$, $n\geq 12$. As a consequence, the existence problem of integer relative $\H_k(n;k)$ is completely settled.
\end{abstract}	
	\maketitle

\section{Introduction and Notation}
Heffter arrays are a class of combinatorial arrays introduced by Archdeacon in \cite{A} which lie at the intersection of multiple areas of Discrete Mathematics, from topological graph theory to graph decompositions and sequencing of groups. Since their introduction, they have attracted considerable attention in research, because of both their wide range of applications and connections with other topics, and due to the challenging problems that arise from them. We refer the interested reader to the survey \cite{DP}.

Relative Heffter arrays have been introduced in \cite{CMPP} as a generalization of the classical concept of Heffter arrays in order to extend their  definition in the context of relative difference families, and further studied in other papers, see for instance
\cite{CPfold,CPPBiembeddings,MP1,MPfold,MT,JMP}. 

Classical Heffter arrays may also be viewed as pairs of orthogonal Heffter
systems. This point of view has recently led to the introduction of Heffter spaces, which are sets of mutually orthogonal Heffter systems \cite{BP}. A relative version of this notion,
generalizing relative Heffter arrays and Heffter spaces, has also been studied
in \cite{JMP}.
 Here we recall the definition:
\begin{definition} \label{def:relative_Heffter}
	Let $t$ be a positive integer dividing $2nk$, and let $J = \langle \frac{2nk+t}{t}\rangle$ be the subgroup of $\Z_{2nk+t}$ of order $t$. 
	A \emph{Heffter array $A$ over $\Z_{2nk+t}$ relative to $J$}, denoted by $\H_t(n; k)$, is an $n\times n$ partially filled  array
	with elements in $\Z_{2nk+t}$ such that:
	\begin{itemize}
		\item[$(\rm{a})$] each row and column contains $k$ filled cells;
		\item[$(\rm{b})$] for every $x \in \Z_{2nk+t}\setminus J$, exactly one between $x$ and $-x$ appears in $A$;
		\item[$(\rm{c})$] every row and column is zero-sum in $\Z_{2nk+t}$.
	\end{itemize}
\end{definition}
A (relative) Heffter array is said to be \textit{integer} if its  rows and columns, viewed as elements in $\pm \{1,\dotsc, \lfloor\frac{2nk+t}{2}\rfloor\}$ (avoiding the multiples of $\frac{2nk+t}{t}$), are zero-sum also in the integers.

\begin{example}\label{example:12}
	An integer relative Heffter array $\H_{5}(12;5)$ over the group $\Z_{125}$, whose subgroup of order $5$ consists of the multiples of $25$:
	\[
	\begin{array}{|c|c|c|c|c|c|c|c|c|c|c|c|}\hline
	-55 & 18 &  &  &  &  &  &  & 43 &  & 7 & -13\\ \hline
	-32 & 52 & 31 &  & -5 &  & -46 &  &  &  &  & \\ \hline
	& -19 & 53 & 17 & -45 & -6 &  &  &  &  &  & \\ \hline
	&  & -33 & 51 & 30 & -44 & -4 &  &  &  &  & \\ \hline
	& -2 & -48 & -20 & 54 & 16 &  &  &  &  &  & \\ \hline
	&  & -3 & -47 & -34 & 62 & 22 &  &  &  &  & \\ \hline
	& -49 &  & -1 &  & -28 & 57 & 21 &  &  &  & \\ \hline
	&  &  &  &  &  & -29 & -56 & 35 & 9 & 41 & \\ \hline
	40 &  &  &  &  &  &  & -15 & -59 & 23 &  & 11\\ \hline
	&  &  &  &  &  &  & 12 & -27 & -60 & 36 & 39\\ \hline
	10 &  &  &  &  &  &  & 38 &  & -14 & -58 & 24\\ \hline
	37 &  &  &  &  &  &  &  & 8 & 42 & -26 & -61\\ \hline
	\end{array}
	\]
\end{example}

In \cite{CMPP} it was shown that an integer relative  Heffter array must satisfy the following necessary conditions:
\begin{proposition} \label{prop:nec_relative}
	Suppose that there exists an integer $\H_t(n;k)$.
	\begin{itemize}
		\item[$\mathrm{(1)}$] If $t$ divides $nk$, then
		\[
		nk \equiv 0 \pmod{4} \text{ or } nk \equiv -t \equiv \pm 1 \pmod{4}.
		\]
		\item[$\mathrm{(2)}$]  If $t = 2nk$, then $k$ must be even.	
		\item[$\mathrm{(3)}$]  If $t  \neq 2nk $ does not divide $nk$, then
		\[
		t + 2nk \equiv 0 \pmod{8}.
		\]
	\end{itemize}
\end{proposition}
In the same article \cite{CMPP}, the authors  show the following existence result:
\begin{theorem}
	Let $3\leq k\leq n$ with $k\neq 5$.
	There exists an integer $\H_k(n;k)$ if and only if one of the following holds:
	\begin{itemize}
		\item[$\mathrm{(1)}$] $k$ is odd and $n\equiv 0,3\pmod 4$;
		\item[$\mathrm{(2)}$] $k\equiv 2\pmod 4$ and $n$ is even;
		\item[$\mathrm{(3)}$] $k\equiv 0\pmod 4$.
	\end{itemize}
	Furthermore, there exists an integer $\H_5(n;5)$ if
	$n\equiv 3\pmod 4$ and it does not exist if $n\equiv 1,2\pmod 4$.
\end{theorem}
Note that  the only case not covered by existence and non-existence results was $k=5$ and $n\equiv 0 \pmod 4$. For this class  only the  two sporadic examples $n=8$ and $n=16$ were given  in \cite{CMPP}.
Our main result is the following.
\begin{theorem}\label{main:thm}
	There exists an integer $\H_5(n;5)$ for every  $n\equiv 0 \pmod 4$, $n \geq 8$.
\end{theorem}
Hence proving that:
\begin{corollary}
Let $3\leq k\leq n$. There exists an integer relative Heffter array
$\H_k(n;k)$ if and only if one of the following holds:
	\begin{itemize}
	\item[$\mathrm{(1)}$] $k$ is odd and $n\equiv0,3\pmod4$;
	\item[$\mathrm{(2)}$] $k\equiv2\pmod4$ and $n$ is even;
	\item[$\mathrm{(3)}$] $k\equiv0\pmod4$.
\end{itemize}
\end{corollary}

In the remainder of this section we fix the notation used throughout the article. 
Given two positive integers $a \leq b$, by $[a,b]$ we mean the set $\{a,a+1,\dotsc, b-1,b\}$.
  For an array $A$ containing elements in $\Z_v$ we represent each entry by the unique integer in $  \left[-\left\lfloor\frac{v-1}{2}\right\rfloor,
  \left\lfloor\frac {v}{2}\right\rfloor\right]$.
   The \textit{support} of $A$, written $\textrm{supp}(A)$, is the set of absolute values of its entries. Following the standard notation present in the literature on Heffter arrays, given an array $A$ and a positive integer $x$, by $A \pm x$ we denote the array where $x$ is added to positive entries and is subtracted from negative entries of $A$. The notation can be naturally applied to the single entries of $A$.

Since many constructions of (relative) Heffter arrays follow a diagonal structure, it is useful to recall a standard notation for the diagonal of the array.
Given a square $n\times n$ array, its diagonals are indexed modulo $n$ by
\[
D_r=\{(i,j): j-i\equiv r-1\pmod n\}.
\]
As an example, the main diagonal is written as $D_1$.

In \cite{DW} the authors introduced the following procedure in order to describe how to fill the cells of an array along its diagonals. This notation has been used in several subsequent papers, see for instance \cite{CMP,CMPP,CPPBiembeddings}.

Let $A$ be an $n \times n$ p.f. array; then, the procedure $diag(r,c,s,\Delta_1,\Delta_2,\ell)$ fills the entries
\[
A(r+i\Delta_1,c+i\Delta_1)=s+i\Delta_2\qquad \textrm{for}\ i\in[0,\ell-1].
\]
The parameters used in the $diag$ procedure have the following meaning:
\begin{itemize}
	\item[-] $r$ denotes the starting row,
	\item[-] $c$ denotes the starting column,
	\item[-] $s$ denotes the entry $A(r,c)$,
	\item[-] $\Delta_1$ denotes the increasing value of the row and column at each step,
	\item[-] $\Delta_2$ denotes how much the entry is changed at each step,
	\item[-] $\ell$ is the length of the chain.
\end{itemize}

\section{Construction of $\H_5(n;5)$ for $n \equiv 0 \pmod{4}, \ n \geq 16$}
This section is devoted to show the construction of an integer relative Heffter array
$\H_{5}(n;5)$ for every $n\equiv 0 \pmod 4$, $n\geq 16$. Combined with the results on the small case $n=8$ and the $\H_5(12;5)$ of Example \ref{example:12}, we obtain the main statement of Theorem \ref{main:thm}. We remark that our proof is inspired by the construction of strictly weak integer Heffter arrays given in \cite{CMP}. 

Since we have to split the proof in two cases, depending on the value of $n \pmod{8}$, in Examples \ref{example:0} and \ref{example:4} we show two step-by-step constructions of $\H_5(n;5)$ for $n=16$ and $n=20$.

We begin by considering the integer $\H_3(n;3)$ with $n\equiv 0\pmod 4$ constructed in Proposition 5.3 of \cite{CMPP},
so let $A$ be the $n \times n$ array built using the following procedures labeled $\texttt{A}$ to $\texttt{J}$:
$$\begin{array}{lcl}
	\texttt{A}:\; diag\left(2,2,1,1,1,\frac{n-4}{2}\right); & \hfill &
	\texttt{B}:\; diag\left(\frac{n+6}{2},\frac{n+6}{2},-\frac{n+4}{2},1, -1,\frac{n-4}{2}\right);\\[3pt]
	\texttt{C}:\; diag\left(2,1,-\frac{5n+4}{2},2,-1,\frac{n}{4}\right); & &
	\texttt{D}:\; diag\left(3,2,-\frac{3n+2}{2},2,-1,\frac{n-4}{4}\right);\\[3pt]
	\texttt{E}:\; diag\left(1,2,\frac{3n}{2},2,-1,\frac{n}{4}\right); &&
	\texttt{F}:\; diag\left(2,3,\frac{5n+2}{2},2,-1,\frac{n-4}{4}\right);\\[3pt]
	\texttt{G}:\; diag\left(\frac{n+6}{2},\frac{n+4}{2},-\frac{5n}{4},2,1,\frac{n}{4}\right); &&
	\texttt{H}:\; diag\left(\frac{n+8}{2},\frac{n+6}{2},-\frac{9n}{4},2,1,\frac{n-4}{4}\right);\\[3pt]
	\texttt{I}:\; diag\left(\frac{n+4}{2},\frac{n+6}{2},\frac{11n+8}{4},2,1,\frac{n}{4}\right); &&
	\texttt{J}:\; diag\left(\frac{n+6}{2},\frac{n+8}{2},\frac{7n+8}{4},2,1,\frac{n-4}{4}\right).
\end{array}$$
We also fill the following cells of $A$:
$$\begin{array}{lclcl}
	A\left(1,1\right)=-\frac{n-2}{2}, & \quad  &
	A\left(\frac{n}{2},\frac{n}{2}\right)=n, & \quad &
	A\left(\frac{n}{2},\frac{n+2}{2}\right)=\frac{7n+4}{4}, \\[3pt]
	A\left(\frac{n+2}{2},\frac{n}{2}\right)=-\frac{9n+4}{4},& &
	A\left(\frac{n+2}{2},\frac{n+2}{2}\right)=\frac{n+2}{2},& &
	A\left(\frac{n+2}{2},\frac{n+4}{2}\right)=\frac{7n}{4}, \\[3pt]
	A\left(\frac{n+4}{2},\frac{n+2}{2}\right)=-\frac{9n+8}{4}, &&
	A\left(\frac{n+4}{2},\frac{n+4}{2}\right)=-\frac{n}{2}.
\end{array}$$	

We remark that the filled cells of $A$ lie on the three diagonals $D_1, D_2, \ D_n$, and 
\[
\textrm{supp}(D_1)=[1,n], \qquad \textrm{supp}(D_n\cup  D_2)=[n+1,3n+1]\setminus\{2n+1\}.
\]
Moreover, the entries of the main diagonal $D_1$ have positive sign in the rows $\{2,3,\dotsc, \frac{n}{2}+1\}$. Let $P$ and $N$  be the following two orderings of the row/column indices of, respectively, the positive and negative entries of $D_1$:
\[
\begin{array}{rlrl}
P&=\left(2,3,\ldots,\frac{n}{2}+1\right) \qquad & 	N&=\left(\frac{n}{2}+3,\frac{n}{2}+4,\ldots,n,1,\frac{n}{2}+2\right), \\
&=\left(p_1,p_2,\dotsc, p_{\frac{n}{2}}\right) &  	&=\left(q_1,q_2,\dotsc, q_{\frac{n}{2}}\right), \\
\end{array}
\]
Let  $E = (e_1,e_2,\dotsc, e_{\frac{n}{2}})$ and $F = (f_1,f_2,\dotsc, f_{\frac{n}{2}})$ respectively be the absolute values of the diagonal entries  of $A$  corresponding to the elements of  $P$ and $N$. Let $\pi$ be the following permutation on the indices $\{1,2,\dotsc, \frac{n}{2}\}$:
\begin{itemize}
	\item if $n \equiv 0 \pmod{8}$:
	\[
	\pi({2r-1}) = {2r}, \qquad \pi({2r}) = 2r-1, \qquad \text{for $r = 1,\dotsc, \frac{n}{8}$}
	\]
	and leaves the remaining elements fixed.
	\item if $n \equiv 4 \pmod{8}$:
	\[
	\begin{aligned}
	&\pi(1)=2, \qquad \pi(2)=3, \qquad \pi(3)=1; \\
	&\pi({2r+2}) = {2r+3}, \qquad \pi({2r+3}) = 2r+2, \qquad \text{for $r = 1,\dotsc, \frac{n-12}{8}$}
	\end{aligned}
	\]
	and leaves the remaining elements fixed.
\end{itemize}
Let now $B$ be the matrix obtained from $A$ by replacing the entries in the main diagonal of $A$ as follows:
\begin{equation} \label{eq_substitution}
\begin{aligned}
	e_i \leftarrow e_{\pi(i)}+4n+2, \qquad 
	f_i\leftarrow -(f_{\pi(i)}+4n+2), \qquad \text{ for $i \in \left[1,\frac{n}{2}\right]$}. \\
\end{aligned}
\end{equation}
Since we have exclusively changed the elements in the main diagonal of $A$, the sum of the $i$-th row is equal to the sum of the $i$-th column. Using the same indices of \eqref{eq_substitution}, it can be seen that their sum is 
\begin{equation} \label{eq_tot_sum}
	\begin{aligned}
&n \equiv 0 \pmod{8}: \left\{\begin{array}{lll}
	\pm(4n+2)& \quad &\text{ if $i\in \left[\frac{n}{4}+1,\frac{n}{2}\right]$,} \\
	\pm(4n+1)& \quad &\text{ if $i \in \left[1,\frac{n}{4}\right]$ and $i$ is even,} \\
	\pm(4n+3)& \quad& \text{ if $i \in \left[1,\frac{n}{4}\right]$ and $i$ is odd.} \\
\end{array}\right.\\
&n \equiv 4 \pmod{8}: \left\{\begin{array}{lll}
	 \pm 4n &\quad &\text{ if $i =3$,} \\
	 \pm (4n+3) &\quad &\text{ if $i \in \{1,2\} $ or $i \in [4,\frac{n}{4}]$ and $i$ even,} \\
	 \pm (4n+1) &\quad &\text{ if $i \in [4,\frac{n}{4}]$ and $i$ odd,} \\
	 \pm (4n+2) &\quad &\text{ if $i\in \left[\frac{n}{4}+1,\frac{n}{2}\right]$.} \\
\end{array}\right.
\end{aligned}
\end{equation}
Now, starting from the array  $B$ we construct a new array $C$ by repeatedly inserting the following auxiliary matrices $Q(a)$ and $T(a)$:
\[
Q(a)=
\begin{array}{|c|c|}\hline
	4n+2-a&a+1 \\ \hline
	a& 4n+1-a\\\hline
\end{array} \qquad T(a)=
\begin{array}{|c|c|c|}\hline
	4n+2-a& a+1& \\ \hline
	&4n+1-a &a+2\\\hline
	a& &4n-a\\\hline
\end{array} \qquad 
\]

Note that the columns of $Q(a)$ sum to $4n+2$, while its rows respectively add to  $4n+3$ and $4n+1$. Regarding $T(a)$, its rows sum to $4n+3$, $4n+3$ and $4n$, while its columns add to $4n+2$.

Assume first that $n \equiv 0 \pmod{8}$. For $r = 1,\dotsc, \frac{n}{8}$ let  $a_1, a_2,b_1,b_2$ be any four distinct elements of $ \{1,3,5,\dotsc, n-1\}$ that have not yet been used in the construction of $C$:
\begin{itemize}
\item add the block $-Q(a_1)$ in the positions
\[
\begin{array}{ll}
	(p_{2r-1},\; p_{\frac{n}{4}+2r-1}),\quad&(p_{2r-1},\; p_{\frac{n}{4}+2r}), \\
	 (p_{2r},\; p_{\frac{n}{4}+2r-1}), \quad & (p_{2r},\; p_{\frac{n}{4}+2r}).
\end{array}
\] 
and the block $-Q(a_2)^t$ in the positions:
\[
\begin{array}{ll}
	( p_{\frac{n}{4}+2r-1},\; p_{2r-1}),\quad&(p_{\frac{n}{4}+2r-1}, \; p_{2r}), \\
	(p_{\frac{n}{4}+2r},\; p_{2r-1} ), \quad & (p_{\frac{n}{4}+2r},\; p_{2r} ).
\end{array}
\] 
\item add the block $Q(b_1)$ in the positions
\[
\begin{array}{ll}
	(q_{2r-1},\; q_{\frac{n}{4}+2r-1}),\quad&(q_{2r-1},\; q_{\frac{n}{4}+2r}), \\
(q_{2r},\; q_{\frac{n}{4}+2r-1}), \quad & (q_{2r},\; q_{\frac{n}{4}+2r}).
\end{array}
\] 
and the block $Q(b_2)^t$ in the positions:
\[
\begin{array}{ll}
	( q_{\frac{n}{4}+2r-1},\; q_{2r-1}),\quad&(q_{\frac{n}{4}+2r-1}, \; q_{2r}), \\
	(q_{\frac{n}{4}+2r},\; q_{2r-1} ), \quad & (q_{\frac{n}{4}+2r},\; q_{2r} ).
\end{array}
\] 

\end{itemize}
It is easy to see that any considered cell does not belong to $D_1,D_2$ or $D_n$, since the absolute difference between its row and column index, say $i$ and $j$, is such that $j\not \equiv i,i-1,i+1 \pmod{n}$.
Recall that the columns of $Q(a)$ sum to $4n+2$, while its rows respectively add to  $4n+3$ and $4n+1$. Hence, the added blocks compensate the row and column sums of \eqref{eq_tot_sum}, so every row and column of $C$ has zero-sum in the integers.

Regarding the support of $C$, recall that: $\textrm{supp}(D_n\cup  D_2)=[n+1,3n+1]\setminus\{2n+1\}$. Moreover, 
\[
\textrm{supp}\left(\bigcup_{a \in \{1,3,\dotsc, n-1\}} Q(a)\right) = [1,n] \cup [3n+2,4n+1] 
\]
\[
\textrm{supp}(D_1) = [4n+3,5n+2].
\]
Hence overall the support of $C$ is $[1,5n+2] \setminus\{2n+1,4n+2\}$, so $C$ is an integer $\H_5(n;5)$.

Assume now that $n\equiv 4 \pmod{8}$, $n\geq20$. Insert first:
\begin{itemize}
	\item $-T(1)$ in the cells $\{p_1,p_2,p_3\} \times \{p_{\frac{n}{4}+1},p_{\frac{n}{4}+2},p_{\frac{n}{4}+3}\}$ and $-T(4)^t$ in the cells $\{p_{\frac{n}{4}+1},p_{\frac{n}{4}+2},p_{\frac{n}{4}+3}\}\times\{p_1,p_2,p_3\} $. 
	\item $T(7)$ in the cells $\{q_1,q_2,q_3\} \times \{q_{\frac{n}{4}+1},q_{\frac{n}{4}+2},q_{\frac{n}{4}+3}\}$ and $T(10)^t$ in the cells $\{q_{\frac{n}{4}+1},q_{\frac{n}{4}+2},q_{\frac{n}{4}+3}\}\times\{q_1,q_2,q_3\} $. 
\end{itemize}
Then for $r = 1,\dotsc, \frac{n-12}{8}$ let $a_1,a_2,b_1,b_2$ be any four distinct elements of $\{13,15,\dotsc, n-1\}$ that have not yet been used into $C$:
\begin{itemize}
	\item 
	add the block $-Q(a_1)$ in the positions
	\[
	\begin{array}{ll}
		(p_{2r+2},\; p_{\frac{n}{4}+2r+2}),\;&(p_{2r+2},\; p_{\frac{n}{4}+2r+3}), \\
		(p_{2r+3},\; p_{\frac{n}{4}+2r+2}), \; & (p_{2r+3},\; p_{\frac{n}{4}+2r+3}).
	\end{array}
	\] 
	and the block $-Q(a_2)^t$ in the positions:
	\[
	\begin{array}{ll}
		( p_{\frac{n}{4}+2r+2},\; p_{2r+2}),\quad&(p_{\frac{n}{4}+2r+2}, \; p_{2r+3}), \\
		(p_{\frac{n}{4}+2r+3},\; p_{2r+2} ), \quad & (p_{\frac{n}{4}+2r+3},\; p_{2r+3} ).
	\end{array}
	\] 
	\item	add the block $Q(b_1)$ in the positions
	\[
\begin{array}{ll}
	(q_{2r+2},\; q_{\frac{n}{4}+2r+2}),\quad&(q_{2r+2},\; q_{\frac{n}{4}+2r+3}), \\
	(q_{2r+3},\; q_{\frac{n}{4}+2r+2}), \quad & (q_{2r+3},\; q_{\frac{n}{4}+2r+3}).
\end{array}
\] 
	and the block $Q(b_2)^t$ in the positions:
	\[
\begin{array}{ll}
	( q_{\frac{n}{4}+2r+2},\; q_{2r+2}),\quad&(q_{\frac{n}{4}+2r+2}, \; q_{2r+3}), \\
	(q_{\frac{n}{4}+2r+3},\; q_{2r+2} ), \quad & (q_{\frac{n}{4}+2r+3},\; q_{2r+3} ).
\end{array}
\] 
	
\end{itemize}
As in the case $n \equiv 0 \pmod{8}$, the cells mentioned above do not belong to $D_1$, $D_2$ or $D_n$, hence the array $C$ has $5$ filled cells for each row and each column. Moreover, the added blocks compensate the row and column sums of \eqref{eq_tot_sum}, hence every row and column of $C$ has zero-sum in the integers.

We conclude by examining the support of $C$: as before,  $\textrm{supp}(D_n\cup  D_2)=[n+1,3n+1]\setminus\{2n+1\}$. On the other hand, 
\[
\textrm{supp}\left(\bigcup_{a \in \{1,4,7,10\}}T(a) \cup \bigcup_{a \in\{13,15,\dotsc, n-1\}} Q(a)\right) = [1,n] \cup [3n+2,4n+1] 
\]
\[
\textrm{supp}(D_1) = [4n+3,5n+2].
\]
So, $C$ is an integer $\H_5(n;5)$.

\begin{example}\label{example:0}
	In this example, we build an $\H_5(n;5)$ for $n=16$ by following the previous construction. We begin with a $\H_3(16;3)$, as described at the beginning of the section:
	\[
	\begin{scriptsize}
	\begin{array}{|c|c|c|c|c|c|c|c|c|c|c|c|c|c|c|c|} \hline
	-7 & 24 &  &  &  &  &  &  &  &  &  &  &  &  &  & -17\\ \hline
	-42 & 1 & 41 &  &  &  &  &  &  &  &  &  &  &  &  & \\ \hline
	& -25 & 2 & 23 &  &  &  &  &  &  &  &  &  &  &  & \\ \hline
	&  & -43 & 3 & 40 &  &  &  &  &  &  &  &  &  &  & \\ \hline
	&  &  & -26 & 4 & 22 &  &  &  &  &  &  &  &  &  & \\ \hline
	&  &  &  & -44 & 5 & 39 &  &  &  &  &  &  &  &  & \\ \hline
	&  &  &  &  & -27 & 6 & 21 &  &  &  &  &  &  &  & \\ \hline
	&  &  &  &  &  & -45 & 16 & 29 &  &  &  &  &  &  & \\ \hline
	&  &  &  &  &  &  & -37 & 9 & 28 &  &  &  &  &  & \\ \hline
	&  &  &  &  &  &  &  & -38 & -8 & 46 &  &  &  &  & \\ \hline
	&  &  &  &  &  &  &  &  & -20 & -10 & 30 &  &  &  & \\ \hline
	&  &  &  &  &  &  &  &  &  & -36 & -11 & 47 &  &  & \\ \hline
	&  &  &  &  &  &  &  &  &  &  & -19 & -12 & 31 &  & \\ \hline
	&  &  &  &  &  &  &  &  &  &  &  & -35 & -13 & 48 & \\ \hline
	&  &  &  &  &  &  &  &  &  &  &  &  & -18 & -14 & 32\\ \hline
	49 &  &  &  &  &  &  &  &  &  &  &  &  &  & -34 & -15\\ \hline
	\end{array}
	\end{scriptsize}
	\] 
	We have:
	\[
	\begin{array}{ll}
		P = (2,3,4,5,6,7,8,9) = (p_1,\dotsc, p_8)  &N = (11,12,13,14,15,16,1,10) = (q_1,\dotsc, q_8)\\
		E = (1,2,3,4,5,6,16,9)= (e_1,\dotsc, e_8)  &F = (10,11,12,13,14,15,7,8)= (f_1,\dotsc, f_8)
	\end{array}
	\]
Let $\pi$ be the following permutation on $[1,8]$:\[
\pi  = (1 \ 2) (3 \ 4)(5)(6)(7)(8).\] 
We construct then the array $B$ by replacing the entries of the main diagonal of $A$ as follows:
	\[
	\begin{array}{ll}
		e_1 \leftarrow e_{\pi(1)}+66=68, \qquad &
		f_1\leftarrow -(f_{\pi(1)}+66)=-77, \\
		e_2 \leftarrow e_{\pi(2)}+66=67,& 
		f_2\leftarrow -(f_{\pi(2)}+66)=-76, \\
		e_3 \leftarrow e_{\pi(3)}+66=70, & 
		f_3\leftarrow -(f_{\pi(3)}+66)=-79, \\
		e_4 \leftarrow e_{\pi(4)}+66=69, &
		f_4\leftarrow -(f_{\pi(4)}+66)=-78, \\
		e_5 \leftarrow e_{\pi(5)}+66=71, & 
		f_5\leftarrow -(f_{\pi(5)}+66)=-80, \\
		e_6 \leftarrow e_{\pi(6)}+66=72, & 
		f_6\leftarrow -(f_{\pi(6)}+66)=-81, \\
		e_7 \leftarrow e_{\pi(7)}+66=82,&
		f_7\leftarrow -(f_{\pi(7)}+66)=-73, \\
		e_8 \leftarrow e_{\pi(8)}+66=75, &
		f_8\leftarrow -(f_{\pi(8)}+66)=-74, \\
	\end{array}
	\]
	obtaining the following:
		\[
	\begin{scriptsize}
		\begin{array}{|c|c|c|c|c|c|c|c|c|c|c|c|c|c|c|c|} \hline
-73 & 24 &  &  &  &  &  &  &  &  &  &  &  &  &  & -17\\ \hline
-42 & 68 & 41 &  &  &  &  &  &  &  &  &  &  &  &  & \\ \hline
& -25 & 67 & 23 &  &  &  &  &  &  &  &  &  &  &  & \\ \hline
&  & -43 & 70 & 40 &  &  &  &  &  &  &  &  &  &  & \\ \hline
&  &  & -26 & 69 & 22 &  &  &  &  &  &  &  &  &  & \\ \hline
&  &  &  & -44 & 71 & 39 &  &  &  &  &  &  &  &  & \\ \hline
&  &  &  &  & -27 & 72 & 21 &  &  &  &  &  &  &  & \\ \hline
&  &  &  &  &  & -45 & 82 & 29 &  &  &  &  &  &  & \\ \hline
&  &  &  &  &  &  & -37 & 75 & 28 &  &  &  &  &  & \\ \hline
&  &  &  &  &  &  &  & -38 & -74 & 46 &  &  &  &  & \\ \hline
&  &  &  &  &  &  &  &  & -20 & -77 & 30 &  &  &  & \\ \hline
&  &  &  &  &  &  &  &  &  & -36 & -76 & 47 &  &  & \\ \hline
&  &  &  &  &  &  &  &  &  &  & -19 & -79 & 31 &  & \\ \hline
&  &  &  &  &  &  &  &  &  &  &  & -35 & -78 & 48 & \\ \hline
&  &  &  &  &  &  &  &  &  &  &  &  & -18 & -80 & 32\\ \hline
49 &  &  &  &  &  &  &  &  &  &  &  &  &  & -34 & -81\\ \hline
		\end{array}
	\end{scriptsize}
	\] 
	Note that the row and column sums are:
	\[
	(-66,  67,  65,  67,  65,  66,  66,  66,  66, -66, -67, -65, -67, -65, -66, -66)\]
	We construct a new array $C$ by inserting into $B$ the following blocks:
	\[
    \begin{array}{rrrr}
	Q(1)=\begin{array}{|c|c|}\hline
		65& 2\\ \hline
		1& 64\\ \hline
	\end{array} & 
	Q(3)=\begin{array}{|c|c|}\hline
		63& 4\\ \hline
		3&62 \\ \hline
	\end{array} &
	Q(5)=\begin{array}{|c|c|}\hline
		61& 6\\ \hline
		5&60 \\ \hline
	\end{array} &
	Q(7)=\begin{array}{|c|c|}\hline
		59&8 \\ \hline
		7&58 \\ \hline
	\end{array}  
	\\
    \\
	Q(9)=\begin{array}{|c|c|}\hline
		57&10 \\ \hline
		9&56 \\ \hline
	\end{array} &
	Q(11)=\begin{array}{|c|c|}\hline
		55&12 \\ \hline
		11&54 \\ \hline
	\end{array} &
	Q(13)=\begin{array}{|c|c|}\hline
		53&14 \\ \hline
		13&52 \\ \hline
	\end{array}& 
	Q(15)=\begin{array}{|c|c|}\hline
		51&16 \\ \hline
		15&50 \\ \hline
	\end{array}  
    \end{array}
	\] 

Specifically, we add the blocks

\begin{center}
\begin{tabular}{r|l}
Block& Cells \\ \hline
    $-Q(1)$ & $\{p_1,p_2\}\times \{p_5,p_6\}=\{(2,6),(2,7),(3,6),(3,7)\}$ \\
     $-Q(3)^t$&  $\{p_5,p_6\} \times \{p_1,p_2\}=\{(6,2),(7,2),(6,3),(7,3)\}$ \\
     $Q(5)$& $\{q_1,q_2\}\times \{q_5,q_6\}=\{(11,15),(11,16),(12,15),(12,16)\}$\\
     $Q(7)^t$&$\{q_5,q_6\} \times \{q_1,q_2\}=\{(15,11),(16,11),(15,12),(16,12)\}$ \\
     $-Q(9)$& $\{p_3,p_4\}\times \{p_7,p_8\}=\{(4,8),(4,9),(5,8),(5,9)\}$\\
     $-Q(11)^t$&$\{p_7,p_8\}\times \{p_3,p_4\}=\{(8,4),(9,4),(8,5),(9,5)\}$ \\
     $Q(13)$&$\{q_3,q_4\}\times \{q_7,q_8\}=\{(13,1),(13,10),(14,1),(14,10)\}$ \\
     $Q(15)^t$&$\{q_7,q_8\}\times \{q_3,q_4\}=\{(1,13),(10,13),(1,14),(10,14)\}$ \\
\end{tabular}
\end{center}

Overall, we construct the following array $C$:
		\[
\begin{scriptsize}
	\begin{array}{|c|c|c|c|c|c|c|c|c|c|c|c|c|c|c|c|} \hline
-73 & 24 &  &  &  &  &  &  &  &  &  &  & 51 & 15 &  & -17\\ \hline
-42 & 68 & 41 &  &  & -65 & -2 &  &  &  &  &  &  &  &  & \\ \hline
& -25 & 67 & 23 &  & -1 & -64 &  &  &  &  &  &  &  &  & \\ \hline
&  & -43 & 70 & 40 &  &  & -57 & -10 &  &  &  &  &  &  & \\ \hline
&  &  & -26 & 69 & 22 &  & -9 & -56 &  &  &  &  &  &  & \\ \hline
& -63 & -3 &  & -44 & 71 & 39 &  &  &  &  &  &  &  &  & \\ \hline
& -4 & -62 &  &  & -27 & 72 & 21 &  &  &  &  &  &  &  & \\ \hline
&  &  & -55 & -11 &  & -45 & 82 & 29 &  &  &  &  &  &  & \\ \hline
&  &  & -12 & -54 &  &  & -37 & 75 & 28 &  &  &  &  &  & \\ \hline
&  &  &  &  &  &  &  & -38 & -74 & 46 &  & 16 & 50 &  & \\ \hline
&  &  &  &  &  &  &  &  & -20 & -77 & 30 &  &  & 61 & 6\\ \hline
&  &  &  &  &  &  &  &  &  & -36 & -76 & 47 &  & 5 & 60\\ \hline
53 &  &  &  &  &  &  &  &  & 14 &  & -19 & -79 & 31 &  & \\ \hline
13 &  &  &  &  &  &  &  &  & 52 &  &  & -35 & -78 & 48 & \\ \hline
&  &  &  &  &  &  &  &  &  & 59 & 7 &  & -18 & -80 & 32\\ \hline
49 &  &  &  &  &  &  &  &  &  & 8 & 58 &  &  & -34 & -81\\ \hline
	\end{array}
\end{scriptsize}
\] 
It can be easily verified that $C$ is an integer $\H_5(n;5)$.

\end{example}

\begin{example}\label{example:4}
In this example, we build an $\H_5(20;5)$	 by following the previous construction. We begin with the following $\H_3(20;3)$:

		\[
		\hspace{-50pt}
\begin{scriptsize}
	\begin{array}{|c|c|c|c|c|c|c|c|c|c|c|c|c|c|c|c|c|c|c|c|} \hline
-9 & 30 &  &  &  &  &  &  &  &  &  &  &  &  &  &  &  &  &  & -21\\ \hline
-52 & 1 & 51 &  &  &  &  &  &  &  &  &  &  &  &  &  &  &  &  & \\ \hline
& -31 & 2 & 29 &  &  &  &  &  &  &  &  &  &  &  &  &  &  &  & \\ \hline
&  & -53 & 3 & 50 &  &  &  &  &  &  &  &  &  &  &  &  &  &  & \\ \hline
&  &  & -32 & 4 & 28 &  &  &  &  &  &  &  &  &  &  &  &  &  & \\ \hline
&  &  &  & -54 & 5 & 49 &  &  &  &  &  &  &  &  &  &  &  &  & \\ \hline
&  &  &  &  & -33 & 6 & 27 &  &  &  &  &  &  &  &  &  &  &  & \\ \hline
&  &  &  &  &  & -55 & 7 & 48 &  &  &  &  &  &  &  &  &  &  & \\ \hline
&  &  &  &  &  &  & -34 & 8 & 26 &  &  &  &  &  &  &  &  &  & \\ \hline
&  &  &  &  &  &  &  & -56 & 20 & 36 &  &  &  &  &  &  &  &  & \\ \hline
&  &  &  &  &  &  &  &  & -46 & 11 & 35 &  &  &  &  &  &  &  & \\ \hline
&  &  &  &  &  &  &  &  &  & -47 & -10 & 57 &  &  &  &  &  &  & \\ \hline
&  &  &  &  &  &  &  &  &  &  & -25 & -12 & 37 &  &  &  &  &  & \\ \hline
&  &  &  &  &  &  &  &  &  &  &  & -45 & -13 & 58 &  &  &  &  & \\ \hline
&  &  &  &  &  &  &  &  &  &  &  &  & -24 & -14 & 38 &  &  &  & \\ \hline
&  &  &  &  &  &  &  &  &  &  &  &  &  & -44 & -15 & 59 &  &  & \\ \hline
&  &  &  &  &  &  &  &  &  &  &  &  &  &  & -23 & -16 & 39 &  & \\ \hline
&  &  &  &  &  &  &  &  &  &  &  &  &  &  &  & -43 & -17 & 60 & \\ \hline
&  &  &  &  &  &  &  &  &  &  &  &  &  &  &  &  & -22 & -18 & 40\\ \hline
61 &  &  &  &  &  &  &  &  &  &  &  &  &  &  &  &  &  & -42 & -19\\ \hline
	\end{array}
\end{scriptsize}
\] 
	We have:
\[
\begin{aligned}
	&P = (2,3,4,5,6,7,8,9,10,11) = (p_1,\dotsc, p_{10}) \qquad N=(13,14,15,16,17,18,19,20,1,12) = (q_1,\dotsc, q_{10})\\
	&E = (1,2,3,4,5,6,7,8,20,11)= (e_1,\dotsc, e_{10}) \qquad F = (12,13,14,15,16,17,18,19,9,10)= (f_1,\dotsc, f_{10})
\end{aligned}
\]
We have the following permutation $\pi $ on $[1,10]$:
\[
\pi = (1\ 2 \ 3) (4 \ 5)(6)(7)(8)(9)(10).
\]
We construct then the array $B$ by replacing the entries of the main diagonal of $A$ as follows:
\[
\begin{array}{ll}
	e_1 \leftarrow e_{\pi(1)}+82=84, \qquad& 
	f_1\leftarrow -(f_{\pi(1)}+82)=-95, \\
	e_2 \leftarrow e_{\pi(2)}+82=85,&
	f_2\leftarrow -(f_{\pi(2)}+82)=-96, \\
	e_3 \leftarrow e_{\pi(3)}+82=83, &
	f_3\leftarrow -(f_{\pi(3)}+82)=-94, \\
	e_4 \leftarrow e_{\pi(4)}+82=87, & 
	f_4\leftarrow -(f_{\pi(4)}+82)=-98, \\
	e_5 \leftarrow e_{\pi(5)}+82=86, &
	f_5\leftarrow -(f_{\pi(5)}+82)=-97, \\
	e_6 \leftarrow e_{\pi(6)}+82=88, &
	f_6\leftarrow -(f_{\pi(6)}+82)=-99, \\
	e_7 \leftarrow e_{\pi(7)}+82=89, & 
	f_7\leftarrow -(f_{\pi(7)}+82)=-100, \\
	e_8 \leftarrow e_{\pi(8)}+82=90, & 
	f_8\leftarrow -(f_{\pi(8)}+82)=-101, \\
		e_9 \leftarrow e_{\pi(9)}+82=102, & 
	f_9\leftarrow -(f_{\pi(9)}+82)=-91, \\
		e_{10} \leftarrow e_{\pi(10)}+82=93, &
	f_{10}\leftarrow -(f_{\pi(10)}+82)=-92, \\
\end{array}
\]
obtaining the following array $B$:
\[
		\hspace{-50pt}
\begin{scriptsize}
	\begin{array}{|c|c|c|c|c|c|c|c|c|c|c|c|c|c|c|c|c|c|c|c|} \hline
-91 & 30 &  &  &  &  &  &  &  &  &  &  &  &  &  &  &  &  &  & -21\\ \hline
-52 & 84 & 51 &  &  &  &  &  &  &  &  &  &  &  &  &  &  &  &  & \\ \hline
& -31 & 85 & 29 &  &  &  &  &  &  &  &  &  &  &  &  &  &  &  & \\ \hline
&  & -53 & 83 & 50 &  &  &  &  &  &  &  &  &  &  &  &  &  &  & \\ \hline
&  &  & -32 & 87 & 28 &  &  &  &  &  &  &  &  &  &  &  &  &  & \\ \hline
&  &  &  & -54 & 86 & 49 &  &  &  &  &  &  &  &  &  &  &  &  & \\ \hline
&  &  &  &  & -33 & 88 & 27 &  &  &  &  &  &  &  &  &  &  &  & \\ \hline
&  &  &  &  &  & -55 & 89 & 48 &  &  &  &  &  &  &  &  &  &  & \\ \hline
&  &  &  &  &  &  & -34 & 90 & 26 &  &  &  &  &  &  &  &  &  & \\ \hline
&  &  &  &  &  &  &  & -56 & 102 & 36 &  &  &  &  &  &  &  &  & \\ \hline
&  &  &  &  &  &  &  &  & -46 & 93 & 35 &  &  &  &  &  &  &  & \\ \hline
&  &  &  &  &  &  &  &  &  & -47 & -92 & 57 &  &  &  &  &  &  & \\ \hline
&  &  &  &  &  &  &  &  &  &  & -25 & -95 & 37 &  &  &  &  &  & \\ \hline
&  &  &  &  &  &  &  &  &  &  &  & -45 & -96 & 58 &  &  &  &  & \\ \hline
&  &  &  &  &  &  &  &  &  &  &  &  & -24 & -94 & 38 &  &  &  & \\ \hline
&  &  &  &  &  &  &  &  &  &  &  &  &  & -44 & -98 & 59 &  &  & \\ \hline
&  &  &  &  &  &  &  &  &  &  &  &  &  &  & -23 & -97 & 39 &  & \\ \hline
&  &  &  &  &  &  &  &  &  &  &  &  &  &  &  & -43 & -99 & 60 & \\ \hline
&  &  &  &  &  &  &  &  &  &  &  &  &  &  &  &  & -22 & -100 & 40\\ \hline
61 &  &  &  &  &  &  &  &  &  &  &  &  &  &  &  &  &  & -42 & -101\\ \hline
	\end{array}
\end{scriptsize}
\] 
The row and column sums of $B$ are:
\[
(-82,   83,   83,   80,   83,   81,   82,   82,   82,   82,   82,  -82,  -83,  -83,  -80,  -83,  -81,  -82,  -82,  -82)
\]
We construct a new array $C$ by inserting into $B$ the following blocks:
\[
	T(1)=\begin{array}{|c|c|c|}\hline
	81& 2&\\ \hline
	& 80&3\\ \hline
	1&&79 \\ \hline
\end{array} \quad 	T(4)=\begin{array}{|c|c|c|}\hline
78& 5&\\ \hline
& 77&6\\ \hline
4&&76 \\ \hline
\end{array} \quad 
	T(7)=\begin{array}{|c|c|c|}\hline
	75& 8&\\ \hline
	& 74&9\\ \hline
	7&&73 \\ \hline
\end{array} \quad 	T(10)=\begin{array}{|c|c|c|}\hline
	72& 11&\\ \hline
	& 71&12\\ \hline
	10&&70 \\ \hline
\end{array}  
\]
and
\[
	Q(13)=\begin{array}{|c|c|}\hline
	69& 14\\ \hline
	13& 68\\ \hline
\end{array} \qquad 
	Q(15)=\begin{array}{|c|c|}\hline
	67& 16\\ \hline
	15& 66\\ \hline
\end{array} \qquad 
	Q(17)=\begin{array}{|c|c|}\hline
	65& 18\\ \hline
	17& 64\\ \hline
\end{array} \qquad 
	Q(19)=\begin{array}{|c|c|}\hline
	63& 20\\ \hline
	19& 62\\ \hline
\end{array} \qquad 
\]
We insert the blocks as follows:

\begin{center}
\begin{tabular}{r|l}
Block& Cells \\ \hline
$-T(1)$& $\{p_1,p_2,p_3\}\times \{p_6,p_7,p_8\}=\{2,3,4\}\times \{7,8,9\}$\\
$-T(4)^t$& $\{p_6,p_7,p_8\}\times \{p_1,p_2,p_3\}=\{7,8,9\}\times \{2,3,4\}$\\
$T(7)$& $\{q_1,q_2,q_3\}\times \{q_6,q_7,q_8\}=\{13,14,15\}\times \{18,19,20\}$\\
$T(10)^t$&$\{q_6,q_7,q_8\}\times \{q_1,q_2,q_3\}=\{18,19,20\}\times \{13,14,15\}$ \\
$-Q(13)$& $\{p_4,p_5\}\times \{p_9,p_{10}\}=\{5,6\}\times \{10,11\}$ \\
$-Q(15)^t$&$ \{p_9,p_{10}\} \times \{p_4,p_5\}= \{10,11\}\times \{5,6\}$ \\
$Q(17)$&$\{q_4,q_5\}\times \{q_9,q_{10}\}=\{16,17\}\times \{1,12\}$ \\
$Q(19)^t$& $ \{q_9,q_{10}\}\times \{q_4,q_5\}=\{1,12\}\times\{16,17\}$\\
\end{tabular}
\end{center}

obtaining the array, which is an integer $\H_5(20;5)$:
\[
\hspace{-50pt}
\begin{scriptsize}
	\begin{array}{|c|c|c|c|c|c|c|c|c|c|c|c|c|c|c|c|c|c|c|c|} \hline
-91 & 30 &  &  &  &  &  &  &  &  &  &  &  &  &  & 63 & 19 &  &  & -21\\ \hline
-52 & 84 & 51 &  &  &  & -81 & -2 &  &  &  &  &  &  &  &  &  &  &  & \\ \hline
& -31 & 85 & 29 &  &  &  & -80 & -3 &  &  &  &  &  &  &  &  &  &  & \\ \hline
&  & -53 & 83 & 50 &  & -1 &  & -79 &  &  &  &  &  &  &  &  &  &  & \\ \hline
&  &  & -32 & 87 & 28 &  &  &  & -69 & -14 &  &  &  &  &  &  &  &  & \\ \hline
&  &  &  & -54 & 86 & 49 &  &  & -13 & -68 &  &  &  &  &  &  &  &  & \\ \hline
& -78 &  & -4 &  & -33 & 88 & 27 &  &  &  &  &  &  &  &  &  &  &  & \\ \hline
& -5 & -77 &  &  &  & -55 & 89 & 48 &  &  &  &  &  &  &  &  &  &  & \\ \hline
&  & -6 & -76 &  &  &  & -34 & 90 & 26 &  &  &  &  &  &  &  &  &  & \\ \hline
&  &  &  & -67 & -15 &  &  & -56 & 102 & 36 &  &  &  &  &  &  &  &  & \\ \hline
&  &  &  & -16 & -66 &  &  &  & -46 & 93 & 35 &  &  &  &  &  &  &  & \\ \hline
&  &  &  &  &  &  &  &  &  & -47 & -92 & 57 &  &  & 20 & 62 &  &  & \\ \hline
&  &  &  &  &  &  &  &  &  &  & -25 & -95 & 37 &  &  &  & 75& 8 & \\ \hline
&  &  &  &  &  &  &  &  &  &  &  & -45 & -96 & 58 &  &  &  & 74 & 9\\ \hline
&  &  &  &  &  &  &  &  &  &  &  &  & -24 & -94 & 38 &  & 7 &  & 73\\ \hline
65 &  &  &  &  &  &  &  &  &  &  & 18 &  &  & -44 & -98 & 59 &  &  & \\ \hline
17 &  &  &  &  &  &  &  &  &  &  & 64 &  &  &  & -23 & -97 & 39 &  & \\ \hline
&  &  &  &  &  &  &  &  &  &  &  & 72 & & 10 &  & -43 & -99 & 60 & \\ \hline
&  &  &  &  &  &  &  &  &  &  &  & 11 & 71 &  &  &  & -22 & -100 & 40\\ \hline
61 &  &  &  &  &  &  &  &  &  &  &  &  & 12 & 70 &  &  &  & -42 & -101\\ \hline
	\end{array}
\end{scriptsize}
\]
\end{example}

\section*{Acknowledgements}
The author would like to thank Anita Pasotti for useful discussions and suggestions.


\begin{thebibliography}{99}
	\bibitem{A} D.S. Archdeacon,
	\textit{Heffter arrays and biembedding graphs on surfaces},
	Electron. J. Combin. \textbf{22} (2015) \#P1.74.

	\bibitem{BP} M. Buratti and A. Pasotti, \textit{Heffter spaces}, Finite Fields Appl. \textbf{98} (2024), 102464.	
	
	\bibitem{CMP} S. Costa, L. Mella and A. Pasotti, {\textit{Weak Heffter arrays and biembedding graphs on non-orientable surfaces}}, Electron. J. Combin. \textbf{31(1)} (2024), \#P1.8.
	
	\bibitem{CMPP} S. Costa, F. Morini, A. Pasotti and M.A. Pellegrini,
	\textit{A generalization of Heffter arrays}, J. Combin. Des. \textbf{28} (2020), 171--206.

    \bibitem{CPfold}
S. Costa  and A. Pasotti,
\textit{On $\lambda$-fold relative Heffter arrays and biembedding multigraphs on surfaces},
European J. Combin. \textbf{97} (2021), 103370.
	
	\bibitem{CPPBiembeddings} S. Costa, A. Pasotti and M.A. Pellegrini,
	\textit{Relative Heffter arrays and biembeddings},
	Ars Math. Contemp. \textbf{18} (2020), 241--271.
	
	\bibitem{DW} J.H. Dinitz and I.M. Wanless,
	\textit{The existence of square integer Heffter arrays},
	Ars Math. Contemp. \textbf{13} (2017), 81--93.

	\bibitem{JMP} L. Johnson, L. Mella  and A. Pasotti, {\textit{On relative simple Heffter spaces}}, J. Algebraic Combin. \textbf{64} (2026), \#32.
	
	
	\bibitem{MT} L. Mella  and T. Traetta, {\textit{Constructing generalized Heffter arrays via near alternating sign matrices}}, J. Combin. Theory A \textbf{205} (2024), 105873.

\bibitem{MP1} F. Morini  and M.A. Pellegrini,  \emph{On the existence of integer relative Heffter arrays}, Discrete
Math. \textbf{343} (2020), 112088.

\bibitem{MPfold}
F. Morini  and M.A. Pellegrini,
\textit{Magic rectangles, signed magic arrays and integer $\lambda$-fold relative Heffter arrays},
Australas. J. Combin. \textbf{80(2)} (2021), 249--280.
    
	\bibitem{DP} A. Pasotti  and J.H. Dinitz, \textit{A survey of Heffter arrays},  
	Fields Inst. Commun. \textbf{86} (2024), 353--392.   
	
\end{thebibliography}
\end{document}